\documentclass[9pt]{IEEEtran}
\usepackage{cite}
\usepackage{graphicx,color}
\usepackage{multirow}
\usepackage{float}
\usepackage[cmex10]{amsmath}
\usepackage{amssymb}
\usepackage{adjustbox}
\usepackage{subfigure}
\usepackage{bm}
\usepackage{color}
\usepackage{relsize}
\usepackage{enumitem}  

\usepackage[top=2.5cm, bottom=3.0cm, left=2.5cm, right=2.5cm]{geometry}

\everymath{\displaystyle}

\newcommand{\beqnum}{\begin{equation}\begin{array}{lcl}}
\newcommand{\eeqnum}{\end{array}\end{equation}}
\newcommand{\beqnom}{\begin{eqnarray}}
\newcommand{\eeqnom}{\end{eqnarray}}
\newcommand{\beqnc}{\begin{center}\begin{eqnarray}}
\newcommand{\eeqnc}{\end{eqnarray}\end{center}}
\newcommand{\beqnlm}{\begin{equation}\vspace{-.5cm}\begin{array}{lll}}
\newcommand{\eeqnlm}{\end{array}\end{equation}}\vspace{-.5cm}
\newcommand{\beq}{\begin{eqnarray*}}
\newcommand{\eeq}{\end{eqnarray*}}
\newcommand{\bef}{\begin{figure}[tbh!]}
\newcommand{\enf}{\end{figure}}

\newtheorem{montheo}{\bf Theorem}

\newtheorem{madef}{\bf Definition}

\usepackage{graphics} 
\usepackage{epsfig} 

\title{\LARGE \bf
Barrier function-Based Variable Gain Super-Twisting Controller
\thanks{This research was supported by Franche-Comt$\acute{e}$ Regional Council (France) [project RECH-MOB15000008]; PASPA and PAPIIT of UNAM [grant number IN115419] and CONACyT, grant 282013.}
}

\author{Hussein Obeid, Salah Laghrouche,  Leonid Fridman, Yacine Chitour, Mohamed Harmouche

\thanks{H. Obeid and S. Laghrouche are with Femto-ST UMR CNRS, Univ. Bourgogne Franche-Comt\'{e}/UTBM, 90010, Belfort, France (e.mail: hussein.obeid@univ-fcomte.fr; salah.laghrouche@utbm.fr).}

\thanks{L. Fridman is with the Departement of Robotics and Control, Engineering Faculty, Universidad Nacional Aut\'{o}noma de M\'{e}xico (UNAM), D.F 04510. M\'{e}xico (e.mail: leonid.fridman@unam.mx).}
\thanks{Y. Chitour is with the Laboratoire des Signaux et Syst\`{e}mes, Universit\'{e} Paris-Sud, CentraleSup\'{e}lec, Gif-sur-Yvette, Paris, France (e.mail:yacine.chitour@l2s.centralesupelec.fr).}
\thanks{M. Harmouche is with Actility, Paris, France (e.mail:mohamed.harmouche@actility.com).}
}

\begin{document}

\maketitle
\pagenumbering{arabic}

\begin{abstract}
\noindent
In this paper, a variable gain super-twisting algorithm based on a barrier function is proposed for a class of first order disturbed systems with uncertain control coefficient and whose disturbances derivatives are bounded but they are unknown. The specific feature of this algorithm is that it can ensure the convergence of the output variable and maintain it in a predefined neighborhood of zero independent of the upper bound of the disturbances derivatives. Moreover, thanks to the structure of the barrier function, it forces the gain to decrease together with the output variable which yields the non-overestimation of the control gain.\end{abstract}

\section{Introduction}


For systems with matching disturbances, the sliding mode controllers have shown their high efficiency \cite{utkin1992sliding}. Indeed, they provide a
closed-loop insensitivity with respect to bounded matched disturbances and guarantee the finite-time convergence to the sliding surface. However, the discontinuity of sliding mode controllers may cause a undesirable big level of chattering in the systems with fast actuators \cite{perezfridman,boiko2008discontinuous}. 
This major obstacle has been attenuated by some
strategies. For systems with fast actuators, relative degree one and Lipschitz disturbances, the super-twisting controller
\cite{levant1993sliding} is one of the most popular
strategies. 
It allows to achieve a second order sliding mode in  finite-time by using a continuous control signal. 
However, the implementation of the super-twisting controller requires the knowledge of an upper bound of the disturbances derivatives, which is unknown or overestimated in practice. Moreover, even the disturbances derivatives are time varying, it will be desirable to follow their variation.

This problem motivates researchers to develop adaptive sliding mode controllers for the case where
the boundaries of the disturbances exist but they are unknown. The general goal of these techniques is to ensure a dynamical adaptation of the control gains in order to be as small as possible while still sufficient to counteract the disturbances and ensure a sliding mode.

The adaptive sliding mode control approaches which exist in the literature can be broadly split into three classes (\cite{Oliveira2018,doi:10.1080/00207179.2016.1194531}):
\begin{enumerate}[label=(\roman*)]
\item Approaches based on the usage of equivalent control value \cite{utkin2013adaptive,edwards2016adaptive,EDWARDS2016183,7506891,Oliveira2018}.
\item  Approaches based on monotonically increasing the gains \cite{negrete2016second,doi:10.1080/00207179.2015.1116713,5717908}. 
\item Approaches based on increasing and decreasing the gains \cite{plestan2010new,shtessel2012novel,SHTESSEL2017229}. 
\end{enumerate}

Approaches in (i) propose to use the equivalent control as an estimation of the disturbance. The latter consist in increasing the gain to enforce the sliding mode to be reached. Once the sliding mode is achieved, the high frequency control signal is low-pass filtered and used as an information about disturbance in controller gain. The sliding mode controller gain consists in the sum of filtered signal and some constant to compensate possible error between real disturbance and its value estimated by filter. However, the algorithm in \cite{utkin2013adaptive} requires the knowledge of the minimum and maximum allowed values of
the adaptive gain, hence, it requires the information of the upper bound of disturbances derivatives. On the other hand, even the other algorithms \cite{edwards2016adaptive,EDWARDS2016183,7506891,Oliveira2018} do not require theoretically the information of the disturbances derivatives, however, in practice, the usage of low-pass filter requires implicitly the information about this upper bound in order to adequately choose the filter time constant. 

 
Strategies in (ii) consist in increasing the gain until the sliding mode is reached, then the gain is fixed at
this value, ensuring an ideal sliding mode for some interval.
When the disturbance grows, the sliding mode can be lost, so the gain
increases and reach it again. 
This second strategy has two main disadvantages: (a) the gain does not decrease, i.e. it will not follow disturbance when it is decreasing; (b)  one cannot be sure that the sliding mode will never lost because it is not ensured that the disturbance will not
grow anymore.


To overcome the first of these disadvantages, approaches in (iii) have been developed. According to these approaches, the gain increases until the sliding mode is achieved and then decreases until the moment it is lost, i.e. the sliding mode is not reached any more. These approaches ensure the finite-time convergence of the sliding variable to some neighborhood of zero without big overestimation of the gain. The main drawback of these approaches is that the size of the
above mentioned neighborhood and the time of convergence
depend on the unknown upper bound of disturbance which are unknown a-priori.


Recently, novel approaches based on the usage of a monitoring function \cite{Hsu2018} and a barrier function \cite{OBEID2018540} have been proposed to adapt the control gains. However, the first strategy has been only applied for the first order sliding mode controller. Whilst the second one has been applied for both first order sliding mode controller and the Levant's Differentiator \cite{OBEID2017}.  

Inspired by (\cite{OBEID2018540,OBEID2017}), this paper proposes a variable gain super-twisting controller based on a new barrier function. This algorithm can drive the output variable and maintain it in a predefined neighborhood of zero, in the presence of Lipschitz disturbances with unknown boundaries. Compared to the earlier work \cite{OBEID2017}, the class of systems considered in this paper contains an additional uncertainty namely the time-dependent uncertain control coefficient. Moreover, the convergence proof of this algorithm is quite different.   

 The advantages of this suggested algorithm are based on the main features of the barrier functions:

\begin{itemize}
\item The output variable converges in a finite time to a predefined
neighborhood of zero, independently of the bound of the disturbances derivatives, and cannot exceed it.

\item The gain provided by the proposed algorithm is not overestimated. This is due to the reason that the
barrier function forces the gain to decrease together with the output variable.

\item The proposed algorithm does not require neither the upper bound of the
disturbances derivatives nor the use of the low-pass filter.
\end{itemize}
This paper is organized as follows. In Section 2, the problem formulation is
given. Section 3 presents the proposed variable gain super-twisting controller. Finally, some conclusions are drawn in section 4.

The notation $\lfloor x \rceil^{\gamma}$, for $ x,\gamma$ in 
$\mathbb{R}$ with $\gamma\geq 0$, is used to represent $|x|^{\gamma} \,\text{sign}(x)$, where $\text{sign}(x)$ is the set-valued function equal to the sign of $x\neq 0$ and $[-1,1]$  for $x=0$ respectively.

\section{Problem formulation}

Consider the first order system described by

\begin{equation}\label{eq:system}
\dot{s}(t)=\gamma(t)u(t)+\delta_0 (t),  
\end{equation}
where $s(t)\in \mathbb{R}$ is the output variable, $u(t)\in \mathbb{R}$ is the super-twisting controller, and $\gamma$ and $\delta_0$ are Lipschitz disturbances such that, if $\delta:=\frac{\delta_0}{\gamma}$, then one has for $t\geq 0$, 
\begin{equation}\label{eq:dist0}
g\leq \gamma(t)\leq G,\quad |\dot{\delta}(t)|\leq M, 
\end{equation}
 where the constant positive bounds $g,G$ and the upper bound $M$ are unknown.

In the presence of Lipschitz disturbances, the standard super-twisting controller \cite{levant1993sliding} given by
\begin{equation}  \label{eq:second}
\begin{cases}
u(t)= -k_1 \lfloor s(t) \rceil^{1/{2}}  + u_2(t), \\ \dot{u}_2(t) =-k_2 \lfloor s(t) \rceil^{0} ,
\end{cases}
\end{equation}
drives both $s(t)$ and $\dot{s}(t)$ to zero in a finite time, i.e. it provides a second order sliding mode if the control gains $k_1$ and $k_2$ are designed as $k_1=1.5\sqrt{M}$ and $k_2=1.1M$. However, the implementation of this standard super-twisting controller requires the information of the unknwon upper bound $M$. Therefore, to overcome this problem, the following variable gain super-twisting controller is considered  \cite{negrete2016second}

\begin{figure}
\begin{minipage}[b]{1\linewidth}
  \centering

  \centerline{\includegraphics[trim = 7.5cm 2.2cm 2.5cm 1.5cm, clip, width=5 cm]{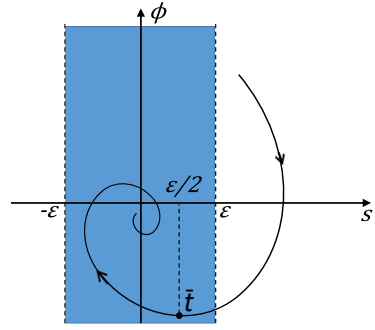}}
\end{minipage}
\vspace*{-0.4cm}
\caption{Phase plane of the proposed super-twisting algorithm}
\label{fig:phase}
\vspace*{1mm}
\end{figure}

\begin{equation}  \label{eq:second}
\begin{cases}
u(t)= - L(t,s_0) \lfloor s(t) \rceil^{1/{2}} + u_2(t), \\ \dot{u}_2(t) =-L^{2}(t,s_0) \lfloor s(t) \rceil^{0},
\end{cases}
\end{equation}
where $s_0:=s(0)$ and $L(t,s_0)$ is the variable gain to be designed in the next section and which depends on the time and  the initial condition $s_0$.
Suppose that $\phi(t)= u_{2}(t)+\delta(t),$ then the dynamic of the first
order system can be expressed as
\begin{equation}  \label{eq:adaptsuper}
\begin{cases}
\dot{s}(t)=\gamma(t)\Big(- L(t,s_0) \lfloor s(t) \rceil^{1/{2}}+\phi(t)\Big), \\ 
\dot{\phi}(t)=- L^{2}(t,s_0) \lfloor s(t) \rceil^{0}+ \dot{\delta}(t).%
\end{cases}%
\end{equation}
The idea behind the proposed algorithm is to first increase the variable gain $L(t,s_0)$ based on the strategy presented in \cite{negrete2016second} until  the output variable reaches the neighborhood of zero $|s|\leq {\varepsilon/{2}}$ at some time $\bar{t}$. Then, for $t>\bar{t}$, the variable gain $L(t,s_0)$ switches to the barrier function and the output variable belongs to the predefined neighborhood of zero $|s(t)|<\varepsilon $.

The trajectory of the proposed super-twisting algorithm in the phase plane $(s(t), \phi (t))$ is illustrated in Fig.~\ref{fig:phase}. It can be shown that, at time $\bar{t}$, the trajectory enters inside the blue vertical strip with a constant bandwidth $]-\varepsilon,\varepsilon[$. The size of this bandwidth remains constant even if the disturbances are time varying. However, when the disturbances change, the size of the vicinity to which converges $s(t)$ changes but the main feature is that it still has an upper bound $\varepsilon $ that it cannot surpass.  


\subsection{Preliminaries }

\subsubsection{Barrier function  }

\begin{madef} \cite{TEE2009918,doi:10.1080/00207179.2011.631192}
Suppose that some $\varepsilon>0$ is given and fixed. For every positive real number $b$, a
barrier function can be defined as an even continuous function  $L_b:x\in (-\varepsilon ,\varepsilon)
\rightarrow L_b(x)\in \left[ b,\infty \right[ $  strictly increasing on  $ \left[ 0,\varepsilon \right[$.
\begin{itemize}

\item $\lim\limits_{|x|\to\varepsilon} L_b(x)= +\infty$.
\item $L_b(x)$ has a unique minimum at zero, i.e. $L_b(0)=b> 0$.
\end{itemize}

\end{madef}

In this paper, the following barrier function is considered
\begin{equation}\label{eq:bf}
L_{b}(x)={\sqrt{\varepsilon} b\over{({\varepsilon- |x|)}^{1/{2}}}},
\quad x\in (-\varepsilon ,\varepsilon),
\end{equation}
where $b$ is a positive constant. 



\section{Main results}

To implement the proposed new algorithm, 
the variable gain $L(t,s_0)$ is defined as follows: first consider the function 
\begin{equation}  
l(t) = L_1 t +L_0 ,  \quad t\geq 0,
\end{equation} 
with $L_0,L_1$ are arbitrary positive constants. Assume first that $|s_0|=|s(0)|\le {\varepsilon/{2}}$. Then, $L(t,s_0)=L_b(s(t))$, with $b=L_0$, as long as the trajectory of \eqref{eq:adaptsuper} is defined. If $|s_0|=|s(0)|> {\varepsilon/{2}}$, then $L(t,s_0)=l(t)$ as long as 
the trajectory $s(\cdot)$ of \eqref{eq:adaptsuper} is defined and 
$|s(t)|>{ \varepsilon /{2}}$. If there exists $\bar{t}(s_0)$ defined as the first time for which $|s(t)|\le{ \varepsilon /{2}}$, then
$L(t,s_0)=L_b(s(t))$, $b=\sqrt{2}l(\bar{t}(s_0))\geq\sqrt{2}L_0$, for $t\geq \bar{t}(s_0)$ and as long as the of trajectory $s(\cdot)$ of \eqref{eq:adaptsuper} is defined.\\
Hence, the variable gain $L(t,s_0)$ is defined, as long as the of trajectory $s(\cdot)$ of \eqref{eq:adaptsuper} is defined, by  
\begin{equation}
\label{eq:adaptgain3}
L(t,s_0)=\begin{cases}
l(t),  \quad \quad \textrm{if} \quad 0\leq t < \bar{t}(s_0),  \\
 L_b(s(t)), \quad \textrm{if}  \quad t\geq\bar{t}(s_0),
\end{cases}
\end{equation} %
with the convention that $\bar{t}(s_0)=0$ if 
$|s(0)|\le \varepsilon/2$ and $\bar{t}(s_0)=\infty$ if $s(\cdot)$ is defined for all times and $|s(\cdot)|> \varepsilon/2$. Since $b=\sqrt{2}l(\bar{t}(s_0))$, then $t\mapsto L(t,s_0)$ defines a continuous function and hence the control $u(\cdot)$ is also continuous as long as it is defined.

\begin{montheo}
	\label{theorem1} 
Let $M$ be the (unknown) upper bound on $|\dot\delta|$, $g,G$ be the bounds on $\gamma(\cdot)$ and $\varepsilon,L_0,L_1>0$ defining the barrier function $L_b$ in \eqref{eq:bf}. Consider System \eqref{eq:adaptsuper} with variable gain $L$ defined in \eqref{eq:adaptgain3}. 
Then, for every $s_0\in\mathbb{R}$ (and initial value of $\dot\delta$), the trajectory of  \eqref{eq:adaptsuper} starting 
at $s_0$ is defined for all non negative times $t$ and there exists a first time $\bar{t}(s_0)$ for which $|s(t)|\le{ \varepsilon/2}$. Then,  
for all $t\ge \bar{t}(s_0)$, one has $|s(t)|<\varepsilon $. Moreover, 
there exists $\nu(M,G,g) >0$ such that, for every trajectory of  \eqref{eq:adaptsuper}, one has that $\limsup_{t\to\infty}|\phi (t)|\leq \nu(M,G,g)$.
\end{montheo}

 The proof of Theorem~\ref{theorem1} is given in Appendix~\ref{proftheor}.

\begin{figure}[tbp]
	\begin{center}
		\subfigure[]{
			\resizebox*{0.95\linewidth}{!}{\label{fig:mcsm}\includegraphics[trim = 0.2cm 6.8cm 0.2cm 7.2cm, clip, width=5cm]{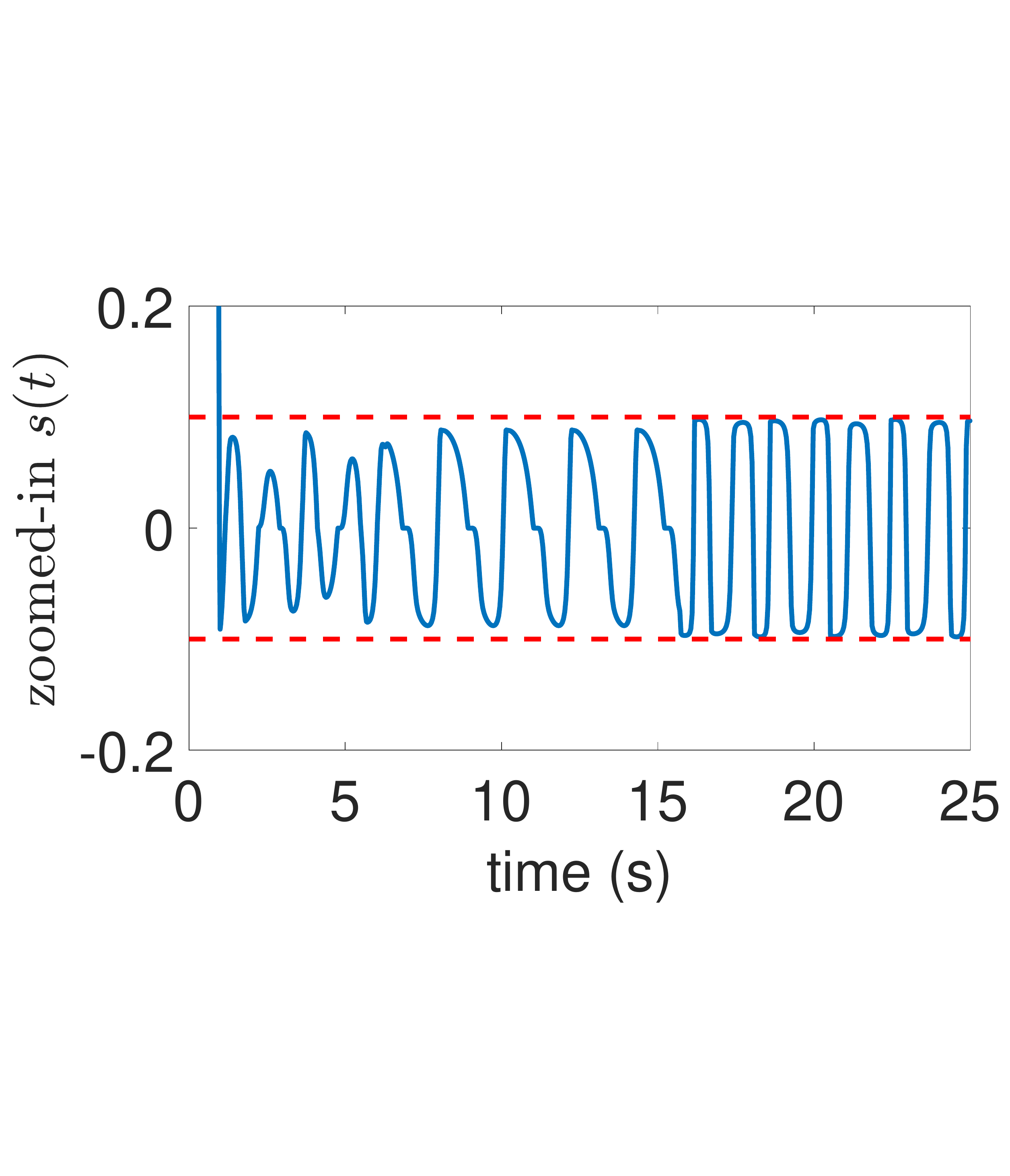}}}\hfill 
		\subfigure[]{
			\resizebox*{0.95\linewidth}{!}{\label{fig:pcsm}\includegraphics[trim = 0.2cm 6.8cm 0.2cm 7.2cm, clip, width=5cm]{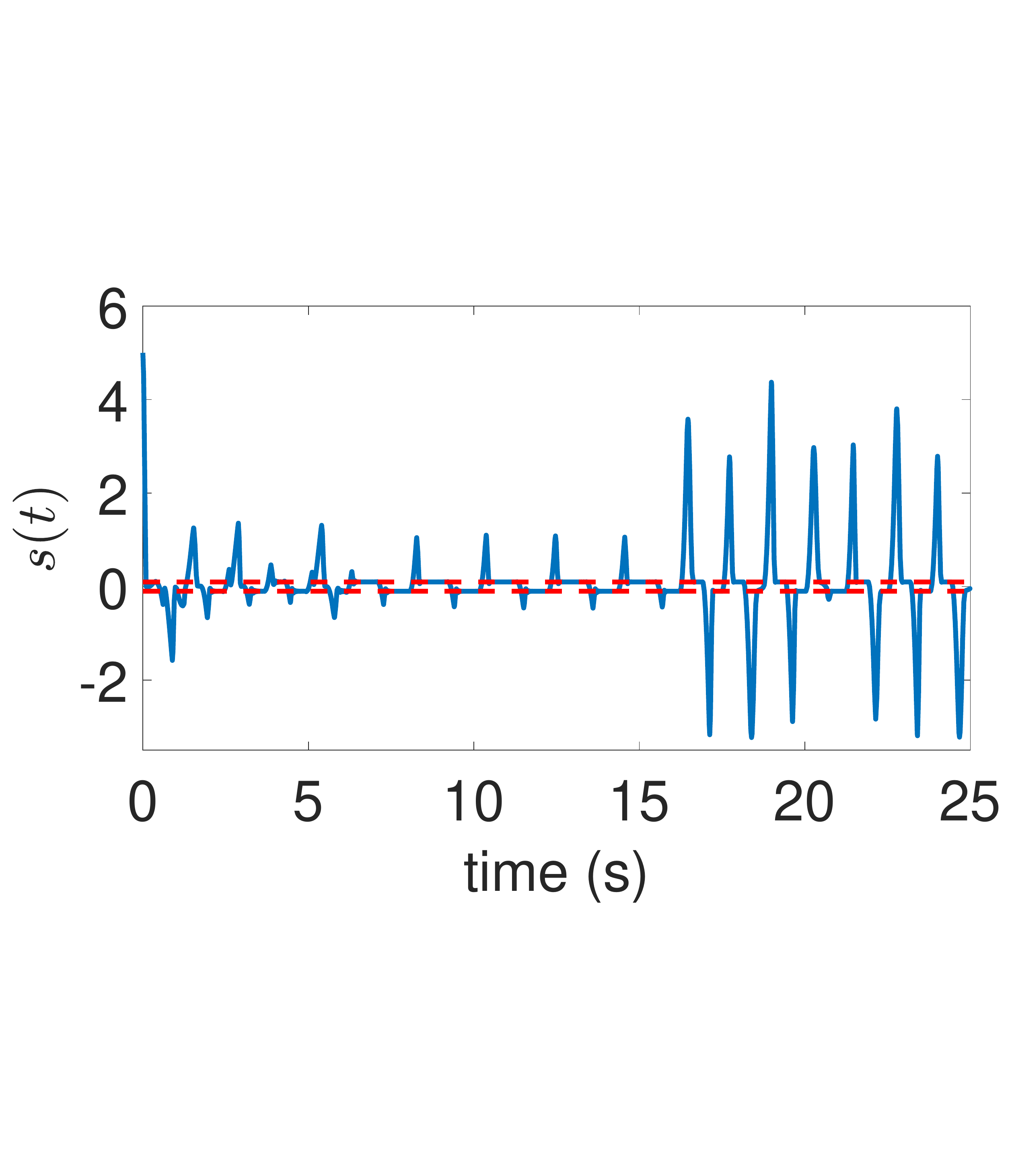}}}%
	\end{center}\vspace*{-0.5cm}
	\caption{Plot of the output variable $s(t)$ with: a) proposed algorithm, b) algorithm presented in \cite{shtessel2012novel}}
	\label{fig:scsm}
\end{figure}

\subsection{Simulation results}

\label{sec:sub5}
The performance of the aforementioned algortihm is compared with the results obtained through the adaptive super-twisting controller presented in \cite{shtessel2012novel}.

In \cite{shtessel2012novel}, the adaptive super-twisting controller is implemented as 
\begin{equation}  \label{eq:shtesselu}
\begin{cases}
u(t)= -\alpha \lfloor s(t) \rceil^{1/{2}} + u_2(t), \\ \dot{u}_2(t) =-{\beta\over{2}} \lfloor s(t) \rceil^{0},
\end{cases}
\end{equation}
where $\beta= 2 \mu \alpha$, and the adaptive gain $\alpha $ is obtained through  
\begin{equation} \label{eq:shtesselalpha}
\dot{\alpha}=    \begin{cases}
     w_1 \sqrt{\gamma_1\over{2}}sign(|s(t)|-\varepsilon), & \text{if}\ \alpha > {\alpha_m} \\     
      \nu  , & \text{if}\ \alpha \le {\alpha_m}\\
     \end{cases}  
  \end{equation}
where $\mu$, $w_1$, $\gamma_1$, $\varepsilon$, $\alpha_m$, and $\nu$ are positive constants to be selected.

In the simulations, $s(0)=5$ and the parameter values of the proposed algorithm are set as $\varepsilon=0.1$, $L_1=1$, $L_0=0.1$. On the other hand, the parameter values of the adaptive super-twisting algorithm \eqref{eq:shtesselu}-\eqref{eq:shtesselalpha} are tuned according to \cite{shtessel2012novel}. Hence, $\mu=1$, $w_1=200$, $\gamma_1=2$, $\nu=\alpha_m=0.01$ while the parameter value of $\varepsilon$ is chosen as the proposed algorithm, i.e. $\varepsilon=0.1$. The disturbances are given by 
\begin{equation}
\gamma(t)= 4+ 2sin(3t).
\end{equation}
\begin{equation}  \label{eq:dist}
	{\delta}(t) = 
	\begin{cases}
		6\sin(5t), \quad \text{if } t \le 2\pi~ s, \\ 
		15\sin(3t),~ \quad \text{if } 2\pi~s < t \le 5\pi~s, \\ 
		30\sin(5t), ~\quad \text{if } t > 5\pi~s.%
	\end{cases}%
\end{equation}
The plots of the output variable $s(t)$ with the proposed algorithm and the algorithm presented in \cite{shtessel2012novel} are compared in Figs.~\ref{fig:mcsm}-\ref{fig:pcsm}. In Figs.~\ref{fig:mcsm} it can be observed that for the proposed algorithm, the output variable does not exceed the predefined neighborhood of zero $\varepsilon =0.1$. On the other hand, it can be noticed in Fig.~\ref{fig:pcsm} that the size of the neighborhood of zero to which converges $s(t)$ with the algorithm presented in \cite{shtessel2012novel} is changing together with the amplitude of disturbances derivatives. Therefore, it cannot be predefined. Moreover, it can be very large when the amplitude of disturbances derivatives is large (for $t> 5\pi~s$, $|s(t)|\le 5= 50~\varepsilon$).

\begin{figure}[tbp]
	\begin{center}
		\subfigure[]{
			\resizebox*{0.95\linewidth}{!}{\label{fig:dosts1}\includegraphics[trim = 0.2cm 6.8cm 0.2cm 7.2cm, clip, width=10cm]{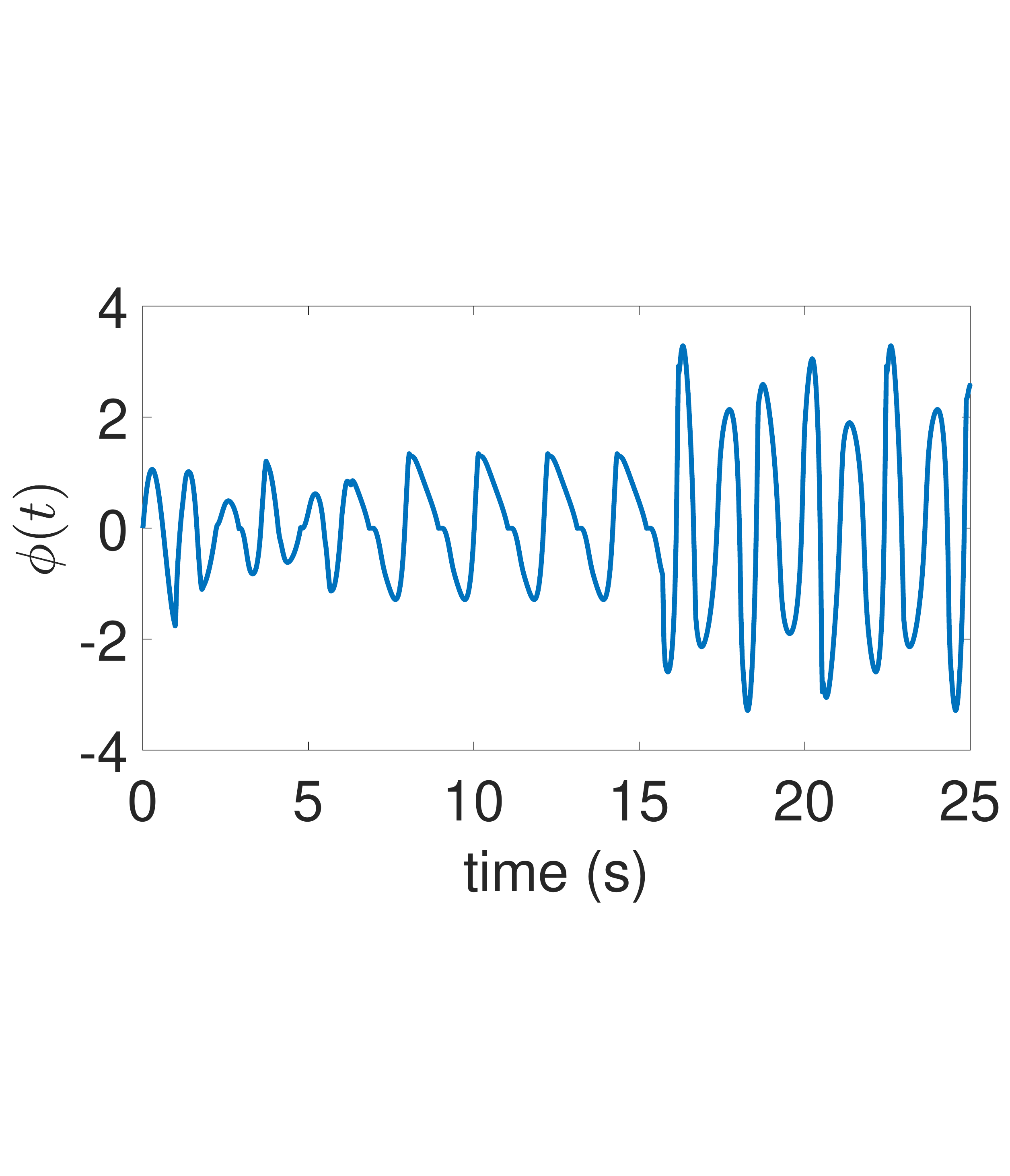}}}%
		\label{fig:barrier2} 
	\hfill 
		\subfigure[]{
			\resizebox*{0.95\linewidth}{!}{\label{fig:dosts3}\includegraphics[trim = 0.2cm 6.8cm 0.2cm 7.2cm, clip, width=10cm]{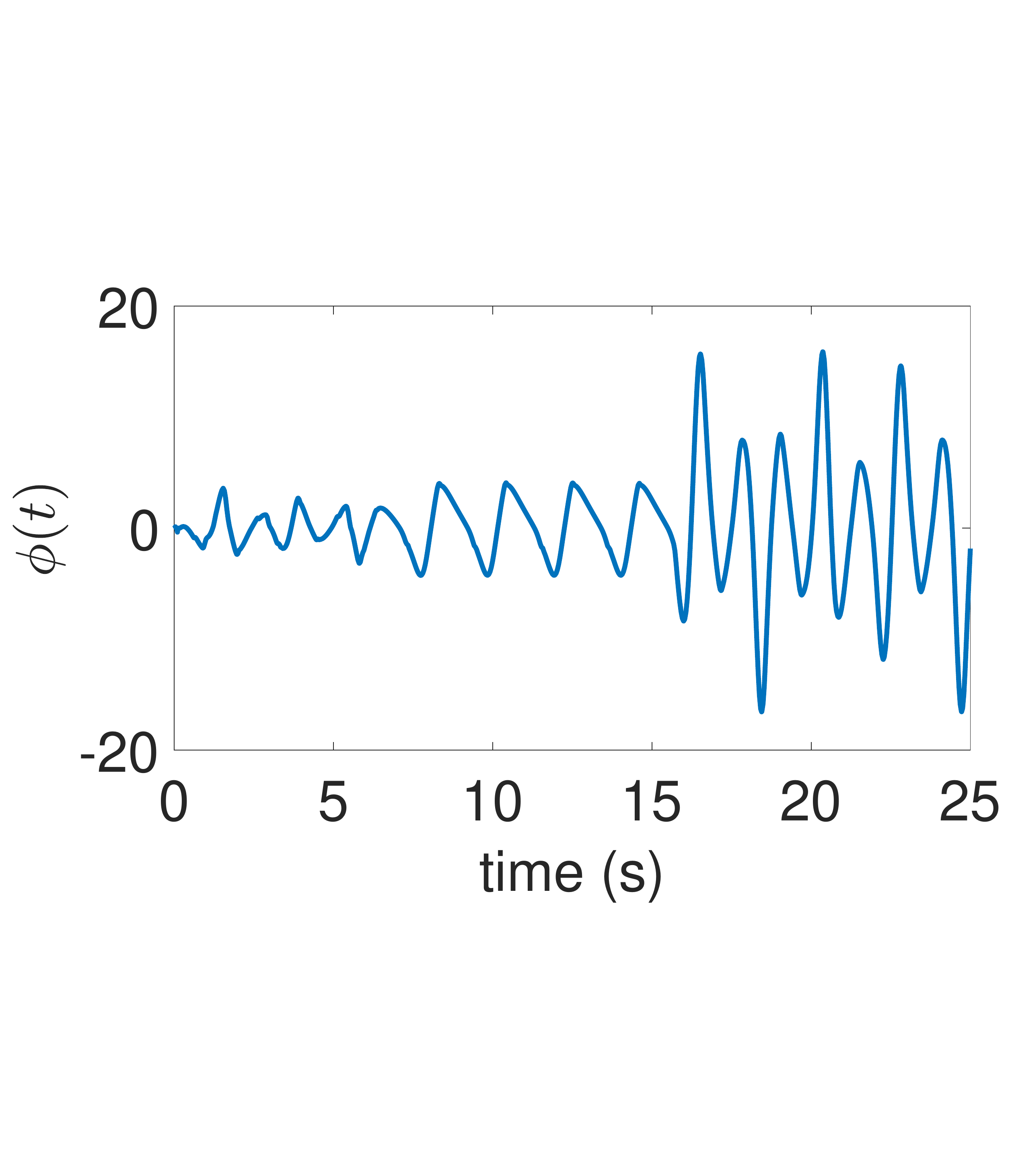}}}%
	\end{center}\vspace*{-0.5cm}
	\caption{Plot of the state variable $\phi(t)$ with: a) proposed algorithm, b) algorithm presented in \cite{shtessel2012novel}  }
	\label{fig:dots}
\end{figure}

\begin{figure}[tbp]
	\begin{center}
		\subfigure[]{
			\resizebox*{0.95\linewidth}{!}{\label{fig:gain1}\includegraphics[trim = 0.2cm 6.8cm 0.2cm 7.2cm, clip, width=5cm]{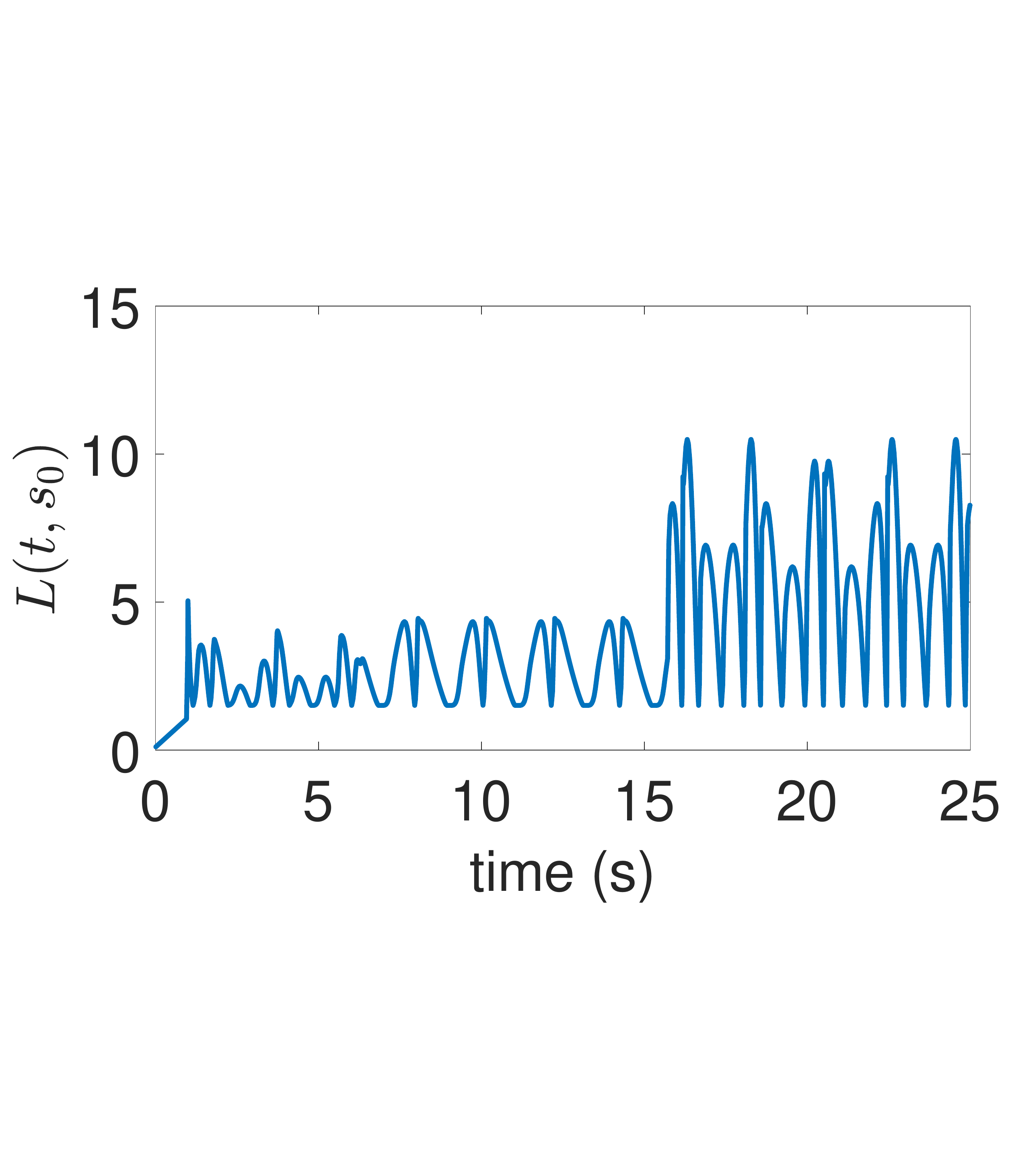}}}\hfill 
		\subfigure[]{
			\resizebox*{0.95\linewidth}{!}{\label{fig:gain2}\includegraphics[trim = 0.2cm 6.8cm 0.2cm 7.2cm, clip, width=5cm]{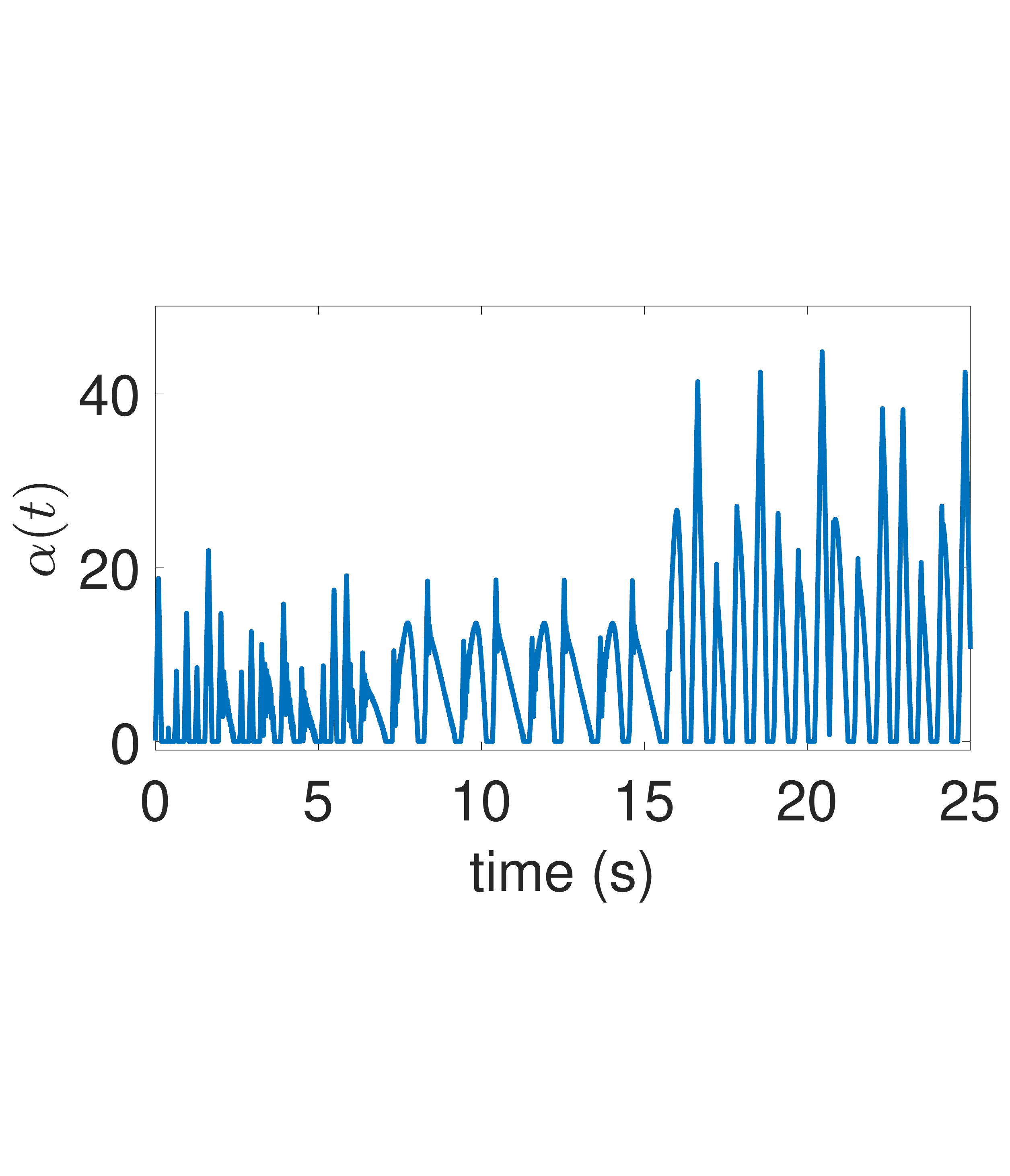}}}%
	\end{center}\vspace*{-0.5cm}
	\caption{Plot of the variable gain with: a) proposed algorithm, b) algorithm presented in \cite{shtessel2012novel}}
	\label{fig:gain}
\end{figure}

Figs.~\ref{fig:dots} illustrates the convergence of the state variable $\phi(t)$ with the proposed algorithm and the algorithm presented in \cite{shtessel2012novel}. It can be confirmed that $\phi (t)$ with the proposed algorithm converges to some vicinity of zero which depends on $M$. Moreover, it can be noticed that the size of vicinity to which converges $\phi(t)$ with the proposed algorithm is less than the one with the algorithm presented in \cite{shtessel2012novel}.

The behaviours of the variable gains with the proposed algorithm and the algorithm presented in \cite{shtessel2012novel} are depicted in Figs.~\ref{fig:gain1}-\ref{fig:gain2}. It can be noticed that both varibale gains can increase and decrease based on the output variable value.

\section{Conclusion}
 
This paper presents a variable gain super-twisting controller for a class of first order disturbed system where the upper bound of the
disturbances derivatives exist, but it is unknown. This algorithm ensures the
convergence of the output variable and prevents its violation outside a predefined neighborhood of zero. Furthermore, the
super-twisting controller gain is not overestimated.

\bibliography{Ref_LHC_TAC}
\begin{appendices}
\section{Proof of Theorem~\ref{theorem1}.} \label{proftheor}
The proof will be done in two steps. 

\textbf{\textit{First step}}:
 It is first shown that there exists a finite time $\bar{t}(s_0)$ for which the output variable $s(t)$, which is part the solution of \eqref{eq:adaptsuper} under the variable gain \eqref{eq:adaptgain3}, becomes $|s(t)|\le {\varepsilon/{2}}$. 
\\
With no loss of generality, we can assume at once that  
$|s(0)| > {\varepsilon/{2}}$. From \eqref{eq:adaptgain3}, the variable gain dynamic is given by $l(t)=L_1t+L_0$ as long as $|s(t)|> {\varepsilon/{2}}$ and the corresponding trajectory of \eqref{eq:adaptsuper} is defined as long as $|s(t)|> {\varepsilon/{2}}$ since, in this case, the growth of the right-hand side of \eqref{eq:adaptsuper} is sublinear with respect to the state variable $(s,\phi)$. Let $I(s_0)$ be the interval where such a dynamics is defined and $I(s_0)$ is of the form $[0,T_0)$. Then, one must prove that $T_0$ is finite, which would at once imply that $T_0$ is the desired time $\bar{t}(s_0)$.\\
Reasoning by contradiction, we assume that $|s|> {\varepsilon/{2}}$ on $I(s_0)$ and we can assume with no loss of generality that $s$ is positive on 
$I(s_0)$. From the second equation of (\ref{eq:adaptsuper}), one gets that 
$$
-l^2(t)-M\leq \dot \phi(t)\leq -l^2(t)+M,\quad\hbox{ on }I(s_0),
$$
which yields easily that $I(s_0)=[0,\infty)$. In particular, we deduce that $\phi(t)$ becomes negative in finite time and remains so. The first equation of (\ref{eq:adaptsuper}) yields that $\dot s(t)\leq -gl(t)\lfloor s(t) \rceil^{1\over{2}}$ and hence convergence of $s$ to zero in finite time. This contradicts the assumption that $|s|> \varepsilon/2$ on $I(s_0)$ and the argument for the first step of the proof of Theorem~\ref{theorem1} is complete.

\textbf{\textit{Second step}}: We now prove that
for all $t \ge \bar{t}(s_0)$, one has $|s(t)|< \varepsilon$ 
and there exists $\nu >0$ only depending on $M,G,g$ such that 
$\limsup_{t\to\infty}|\phi (t)|\leq \nu$. We use $I(s_0)$ to denote the interval of times $t\geq \bar{t}(s_0)$ for which the corresponding trajectory of (\ref{eq:adaptsuper}) is defined. In the argument below, we only work for times $t\in I(s_0)$. In particular, $|s(t)|< \varepsilon$ for $t\in I(s_0)$. It will also be clear for the argument below that $I(s_0)$ is infinite.\\
We consider the following change of variables $y=(y_1,y_2)$ given by
\begin{equation}
y_1 = L_b^2 s(t),\quad
y_2 =\phi(t).
\end{equation}
In theses variables, the system can be written as follows
\begin{equation} \label{eq:systemy}
\begin{cases}
\dot{y}_1=  2 {\dot L_b\over{L_b}} y_1 + L_b^2 \dot{s}, \\
\dot{y}_2= \dot{\phi},
\end{cases}
  \end{equation} 
  and $|y_1(\bar{t}(s_0))|\leq  b^2\varepsilon$.
  Since 
\begin{equation} \label{eq:error}
\dot L_b =  {{1\over{2}} b \sqrt{\varepsilon} {(\varepsilon- |s|)}^{- {1/{2}}} \dot{s} \lfloor s \rceil^{0} \over{(\varepsilon- |s|)}  } ,
  \end{equation} 
then the term ${\dot L_b \over{L_b}}$ can be written as
\begin{equation} \label{eq:dotloverl}
{\dot L_b \over{L_b}}=  {1\over{2}}{ \dot{s} \lfloor s \rceil^{0}  \over{(\varepsilon- |s|)}  } = C_0 L^2_b \dot{s}\lfloor s \rceil^{0} ,
  \end{equation} 
with $C_0=  {1\over{2 b^2 \varepsilon}}$. Substituting \eqref{eq:dotloverl} and \eqref{eq:adaptsuper} into \eqref{eq:systemy}, we obtain 
  \begin{equation} \label{eq:error}
\begin{cases}
{\dot{y}_1\over{(1+  2 C_0 |y_1|)}}=L^2_b\gamma(t)  ( - \lfloor y_1 \rceil^{1/{2}}  + y_2), \\
\dot{y_2}= -  L_b^2 \lfloor y_1 \rceil^{0}  + \dot{\delta}(t) .
\end{cases}
  \end{equation} 
We use a new time scale $\tau$ defined by $\tau(\bar{t}(s_0))=0$ and $d\tau=L_b^2dt$. We use $'$ to denote the derivative with respect to $\tau$.
System \eqref{eq:error} can be rewritten as 
   \begin{equation} \label{eq:error1}
\begin{cases}
{y'_1\over{(1+  2 C_0 |y_1|)}}= \gamma(t) ( -\lfloor y_1 \rceil^{1/{2}} + y_2), \\
y'_2= -  \lfloor y_1 \rceil^{0} + {\dot{\delta}(t) \over{L^2_b}} .
\end{cases}
  \end{equation}   
  Note that 
  $$
  {|\dot{\delta}(t)| \over{L^2_b}}\leq \frac{M}{L_b^2}=\frac{M}{b^2(1+2C_0|y_1|)}\leq\frac{M}{b^2}.
  $$
Then, the second equation of (\ref{eq:error1}) shows that $y'_2$ is bounded and hence $y_2$ has at most a linear growth on $I(s_0)$. Then, from the first equation of (\ref{eq:error1}), one gets that the time derivative of the positive function 
$\ln (1 + 2C_0|y_1|)$ is upper bounded by a function with linear growth on $I(s_0)$, yielding easily that $I(s_0)=[0,\infty)$.\\
  Now consider the following Lyapunov function  
  \begin{equation}\label{eq:v2}
  \begin{split}
V_1(t)&= {\ln (1 + 2C_0|y_1|) \over{2C_0}} \Big(1- {1\over{4}} \lfloor y_1 \rceil^{0} \sigma(y_2) \Big) +{ F(y)y_2^2\over{2}} ,
\end{split}
    \end{equation}  
where $\sigma(y_2)$ is a saturation function defined as 
       \begin{equation} \label{eq:sat}
\sigma(y_2) =\lfloor y_2 \rceil^{0}\min(|y_2|,1),
  \end{equation} 
and $F(y)=g$ if $\lfloor y_1\rceil^{0}y_2\leq 0$ and $G$ if 
$\lfloor y_1\rceil^{0}y_2>0$.
 The following inequality holds  for every $(y_1,y_2)\in\mathbb{R}^2$, 
  \begin{equation*}\label{eq:ineqv2}
 {3 \ln (1 + 2C_0|y_1|) \over{ 8 C_0}}  +{ gy_2^2\over{2}} \le V_1 \le {5\ln (1 + 2C_0|y_1|) \over{8 C_0}} +{ Gy_2^2\over{2}},
    \end{equation*}    
hence $V_1$ is positive-definite and moreover of class $C^1$ {if $y_1\neq 0$  and $\vert y_2\vert\geq 2$.}
The statement of the theorem will be a consequence of the following fact: for every trajectory of \eqref{eq:error1}, one has
\begin{equation}\label{eq:V1}
\limsup_{\tau\to\infty}V_1\leq V^*,
\end{equation}
where $V^*$ is a positive constant only depending on $M,G,g$.\\
  The time derivative of the Lyapunov function \eqref{eq:v2} is 
  given by
   \begin{equation} \label{eq:dotv2}
   \begin{split} 
V'_1 &={y'_1 \lfloor y_1 \rceil^{0} \over{(1 + 2C_0|y_1|)}} +Fy_2y'_2 -  { y'_1 \sigma(y_2) \over{4(1 + 2C_0|y_1|)}} \\&-   { (\sigma(y_2))' \lfloor y_1 \rceil^{0} \ln (1 + 2C_0|y_1|) \over{8 C_0}}.
\end{split} 
    \end{equation}
After easy computations, one gets that     
 \begin{equation} \label{eq:dotv33}
   \begin{split} 
   V'_1 &\leq - {g \over{2}} |y_1|^{\frac12}-|y_2|\Big({g \over{4}} -{M \over{b^2(1+2C_0|y_1|)}}   \Big)\\
&+\sigma'(y_2){M \over{8b^2C_0}}
{\ln (1 + 2C_0|y_1|) \over{1 + 2C_0|y_1|}}\\&+\Big(\gamma(t)-F(y(t))\Big)\lfloor y_1 \rceil^{0}y_2.
\end{split} 
    \end{equation}
Note that a non trivial trajectory of (\ref{eq:error1}) crosses 
the line $y_1=0$ in isolated points. Moreover, at a time $\tau$ where $y(\tau)=(0,y_2(\tau))$ with $y_2(\tau)>0$, $V_1$ admits an isolated discontinuity jumping from $gy^2_2(\tau)$ to $Gy^2_2(\tau)$ and conversely if $y_2(\tau)<0$. Finally, note that, outside of the line $y_1=0$ and for $\vert y_2\vert\geq 2$, $\tau\mapsto V_1(\tau)$ is  a $C^1$ function of the time $\tau$ verifying
 \begin{equation} \label{eq:dotv3}
   \begin{split} 
   V'_1 &\leq - {1 \over{2}} |y_1|^{\frac12}-|y_2|\Big({1 \over{4}} -{M \over{b^2(1+2C_0|y_1|)}}   \Big)\\
&+\sigma'(y_2){M \over{8b^2C_0}}
{\ln (1 + 2C_0|y_1|) \over{1 + 2C_0|y_1|}}.
\end{split} 
    \end{equation}
We will assume without loss of generality that $C_M:=M/b^2$ is much larger than $1$ and hence
$$
y_M:=\frac4{C_0}(C_M-\frac18)>0.
$$
Define 
$$
V^*=8(1+G^2+\frac1{g^2})C^2_M+{5\ln (1 + 2C_0y_M) \over{ C_0}}(1+\frac1g).
$$ 
Let us first prove by contradiction that,
for every trajectory of \eqref{eq:error1}, there exists an increasing sequence of times $(\tau_n)_{n\geq 0}$ tending to infinity such that 
\begin{equation}\label{eq:tau1}
V_1(\tau_n)<\frac{V^*}2,\quad n\geq 0.
\end{equation} 
Hence, if \eqref{eq:tau1} does not hold true for some trajectory $(s(\cdot),\phi(\cdot))$ of \eqref{eq:error1}, there exists $\tau_0\geq 0$ such that one has, for $\tau\geq\tau_0$
$V_1(\tau)\geq V^*/2$.
In the strip $|y_1|\leq y_M$, one has $|y_2|\geq \sqrt{V^*}/2$ and 
$|y'_2+ \lfloor y_1 \rceil^{0}|\leq 1/8$. Moreover, an obvious argument shows that if any trajectory of \eqref{eq:error1} is in the strip at some time $\tau_*\geq \tau_0$
then it must go through it in time $\tau_{cross}$ verifying
\begin{equation}\label{eq:cross}
\tau_{cross}\leq \frac{2\ln(1+4C_M)}{g|y_2(\tau_*)|}\leq
\frac{\ln(1+4C_M)}{C_M}.
\end{equation}
We claim that the trajectory $(s(\cdot),\phi(\cdot))$ of \eqref{eq:error1} must enter into the strip. Indeed, otherwise
$$
V'_1\leq - {1 \over{2}} |y_1|^{\frac12}-{1 \over{8}}|y_2|,
$$
for $\tau$ large enough, i.e. $(V_1^{1/2})'\leq -C_1$ for some positive constant. One would then have convergence in finite time to zero, which is impossible. Hence the trajectory $(s(\cdot),\phi(\cdot))$ of \eqref{eq:error1} must reach a point $(y_M,y_2(\tau_*))$ with $y_2(\tau)\geq \sqrt{V^*}/2$ (the last point due to the symmetry with respect to the origin of trajectories of \eqref{eq:error1}). Then, a simple examination of the phase portrait of \eqref{eq:error1} in the region $y_1>y_M$ shows that 
the trajectory $(s(\cdot),\phi(\cdot))$ must enter the region 
$0<2y_2<y_1^{\frac12}$ and $y_1>y_M$ and exit it in finite time along the axis $y_2=0$. Hence, there exists an interval of time $[\tau_1,\tau_2]$ such that 
$$
\sqrt{y_1}(\tau_1)=2y_2(\tau_1),\quad y_2(\tau_2)=0.
$$
Then $y_2(\tau_1)\geq \sqrt{V^*}/2$, $\tau_2-\tau_1\geq \frac74y_2(\tau_1)$ and 
$$
-\frac{\lfloor y_1 \rceil^{1\over{2}}}2\leq {y'_1\over{(1+  2 C_0 |y_1|)}}\leq -\lfloor y_1 \rceil^{1\over{2}}.
$$
One deduces that $\frac{C_0}2\leq \Big(\frac1{y_1^{1/2}}\Big)'$ and thus 
\begin{equation}\label{eq:tau12}
\frac78C_M\leq \tau_2-\tau_1\leq \frac2{C_0y_1^{1/2}(\tau_2)}\leq 4\exp(-4C_0C_M^2).
\end{equation}
This is clearly impossible and hence \eqref{eq:tau1} is proved.\\
We now prove that \eqref{eq:V1} holds true and the argument goes by contradiction. Indeed, if it were not true, then there exists 
a trajectory $(s(\cdot),\phi(\cdot))$ of \eqref{eq:error1} and an increasing sequence of times $(\tilde{\tau}_n)_{n\geq 0}$ tending to infinity such that 
\begin{equation}\label{eq:tau13}
V_1(\tilde{\tau}_n)>\frac{3V^*}4,\quad n\geq 0.
\end{equation} 
Then, there exists a pair of times (still denoted) $\tau_1<\tau_2$
such that $V_1(\tau_1)=V^*/2$, $V_1(\tau_2)=3V^*/4$ and
\begin{equation}\label{eq:V12}
\frac{V^*}2\leq V_1(\tau)\leq \frac{3V^*}4,\ \tau\in[\tau_1,\tau_2].
\end{equation} 
We will contradict the existence of such a pair of times 
$\tau_1<\tau_2$, which will conclude the proof of Theorem~\ref{theorem1}. Recall that $V_1$ can increase only by going through the strip $|y_1|\leq y_M$. Since the time $\tau_{cross}$ needed to cross that strip is given by \eqref{eq:cross}, one deduces from \eqref{eq:ineqv2} and \eqref{eq:V12} that the increase of $V_1$ by crossing the strip $|y_1|\leq y_M$ is upper bounded by $\ln(1+4C_M)$.
Therefore the trajectory $(s(\cdot),\phi(\cdot))$ of \eqref{eq:error1}
must go through the strip $|y_1|\leq y_M$ at least twice. Combining the above fact with the phase portrait of \eqref{eq:error1} in the region $y_1>y_M$, one deduces that there exists a pair of times 
$\tau'_1<\tau'_2$ in $[\tau_1,\tau_2]$ such that $y_1^{\frac12}(\tau'_1)=2y_2(\tau'_1)$ and $y_2(\tau'_2)=0$. We then arrive at an equation similar to \eqref{eq:tau12} and reach a
contradiction.

\end{appendices}




 

\end{document}